\documentclass[english]{article}
\usepackage[T1]{fontenc}
\usepackage[latin9]{inputenc}
\usepackage{amsmath}
\usepackage{amssymb}
\usepackage{babel}

\begin{document}

\title{Experimental Number Theory\\
Part I : Tower Arithmetic}

\author{by Edinah K. Gnang}

\maketitle

\section{Introduction}

We introduce in this section an Algebraic and Combinatorial approach
to the theory of Numbers. The approach rests on the observation that
numbers can be identified with familiar combinatorial objects namely
\emph{rooted trees}, which we shall here refer to as \emph{towers}.
The bijection between numbers and towers provides some insights into
unexpected connexions between Number theory, combinatorics and discrete
probability theory.

\subsubsection*{Definition 1.1}

Let $\pmb{X}$ denote a $n$ dimensional vector whose entries are
distinct variables defined by \begin{equation}
\pmb{X}=\left(x_{k}\right)_{1\le k\le n}\,,\end{equation}
a \emph{tower expansion }( or simply a tower ) over $\pmb{X}$ is
a finite product of iterated exponentiations over the entries of $\pmb{X}$.
Furthermore the set of towers over the entries of $\pmb{X}$ is denoted
$\mathcal{T}\left(\pmb{X}\right)$.

\subsubsection*{Example 1.1}

Let $\pmb{X}=\left(x,y\right)$, the expressions bellow feature three
towers 

\begin{equation}
x,\quad x^{\left(y^{x}\right)},\qquad x^{\left(y^{x}\right).\left(x^{\left(x^{y}\right)}\right)}.y^{\left(x^{y}\right).\left(y^{\left(y^{y}\right)}\right)}\:.\end{equation}
\\
The \emph{height} of a tower%
\footnote{I shall often omit the parenthesis indicating the order in which the
iterated exponentiations are being carried out, but we shall always
assume that the iterations are performed from the top down .%
} indicates the maximum number of iterated exponentiations occurring
in the tower, the \emph{base} of the tower refers to the bottom level
of the tower. Finally we will call the variables appearing at the
base level of a tower the \emph{pillars} of the tower. \\
\\
\\
\pmb{ \emph{Theorem}: Fundamental Theorem of Arithmetic ( F.
T. A. ):} Every positive integer greater than $1$ can be written
uniquely as a product of powers of primes. ( the expression written
in non decreasing order of the primes.)\\
\\
\pmb{ \emph{Corollary} :} Every positive integer greater than
$1$ can be written uniquely as a tower expansion over the primes
( the primes at each level of the tower are written in increasing
order ).\\
we will not discuss here the proof of the F. T. A. but refer the
reader to a beautiful discussion on the proof of the F. T. A. in \cite{Stein-2008}.

\subsubsection*{Definition 1.2}

A \emph{Formal Tower Series} is a series expansion which consists
of a linear combination of distinct but not necessarily finitely many
towers. The coefficients in the linear combination are assumed to
originate from a field noted here $\mathbb{F}$ ( preferably finite
). The set of such Formal Tower Series is denoted $\mathbb{F}\left[\mathcal{T}\left(\pmb{X}\right)\right]$.
Furthermore a linear combination of finitely many distinct towers
will be called a\textit{ polytower.}

\subsection{Revisiting Euler's product formula for the Riemann zeta function.}

Let $\pmb{X}$ denote an infinite dimensional vector of variables
defined by\begin{equation}
\pmb{X}=\left(x_{k}\right)_{1\le k\le\infty}\end{equation}
As mentioned earlier the F. T. A. induces a bijection between $\mathbb{N}\backslash\left\{ 0,1\right\} $
and towers over the vector of primes, that bijection in turn suggests
a natural bijection between $\mathbb{N}\backslash\left\{ 0,1\right\} $
and $\mathcal{T}\left(\pmb{X}\right)$ as illustrated bellow \begin{equation}
\begin{array}{ccc}
T_{\pmb{X}}\left(2\right) & = & x_{1}\\
T_{\pmb{X}}\left(3\right) & = & x_{2}\\
T_{\pmb{X}}\left(4\right) & = & x_{1}^{x_{1}}\\
T_{\pmb{X}}\left(5\right) & = & x_{3}\\
T_{\pmb{X}}\left(6\right) & = & x_{1}\, x_{2}\\
T_{\pmb{X}}\left(7\right) & = & x_{4}\\
 & \vdots\end{array}\end{equation}
Let us now introduce the binary operator $\pmb{\mathfrak{R}}\left(.\,,\,.\right)$
which we will refer to as the \emph{raiser} \emph{operator} for reason
that will be apparent subsequently.\[
\pmb{\mathfrak{R}}:\:\left\{ T_{\pmb{X}}\left(k\right)\right\} _{1\leq k\leq\infty}\times\mathbb{F}\left[\mathcal{T}\left(\pmb{X}\right)\right]\rightarrowtail\mathbb{F}\left[\mathcal{T}\left(\pmb{X}\right)\right]\]
\begin{equation}
\pmb{\mathfrak{R}}\left(x,\:\sum_{k\in\mathbb{N}/\left\{ 0,1\right\} }a_{k}\: T_{\pmb{X}}\left(k\right)\right)=\sum_{k\in\mathbb{N}\backslash\left\{ 0,1\right\} }a_{k}\: x{}^{T_{\pmb{X}}\left(k\right)}\,.\end{equation}
More generally we write\[
\pmb{\mathfrak{R}}:\left\{ T_{\pmb{X}}\left(k\right)\right\} _{1\leq k\leq\infty}\times\mathbb{F}\left[\mathcal{T}\left(\pmb{X}\right)\right]\rightarrowtail\mathbb{F}\left[\mathcal{T}\left(\pmb{X}\right)\right]\]
\begin{equation}
\pmb{\mathfrak{R}}\left(T_{\pmb{X}}\left(l\right),\:\sum_{k\in\mathbb{N}/\left\{ 0,1\right\} }a_{k}\: T_{\pmb{X}}\left(k\right)\right)=\sum_{k\in\mathbb{N}\backslash\left\{ 0,1\right\} }a_{k}\:\left(T_{\pmb{X}}\left(l\right)\right)^{T_{k}\left(\pmb{X}\right)}.\end{equation}
For most of our discussion however we will require the first of the
two definition of \emph{raiser operator} but we point out that the
first definition follows from the more general following definition.
\\
We recall Euler's product identity for the Riemann Zeta function
as expressed by\begin{equation}
\prod_{k\in\mathbb{N}\backslash\left\{ 0\right\} }\left(1-p_{k}{}^{-s}\right){}^{-1}=1+\sum_{n\in\mathbb{N}\backslash\left\{ 0,1\right\} }n^{-s}\:.\end{equation}
We will come to think of the identity above as expressing an invariance
principle. One important reason for introducing formal tower series
is to validate in some sense identities of the form\begin{equation}
\prod_{k\in\mathbb{N}\backslash\left\{ 0\right\} }\left(\sum_{t\in\mathbb{N}}p_{k}^{t}\right)=1+\sum_{n\in\mathbb{N}\backslash\left\{ 0,1\right\} }n\end{equation}
Here is how the expression above can be thought to be not only meaningful
but also depicting a fundamental invariance principle.\begin{equation}
\prod_{k\in\mathbb{N}\backslash\left\{ 0\right\} }\left(1+\pmb{\mathfrak{R}}\left(x_{k},\;1+\sum_{n\in\mathbb{N}\backslash\left\{ 0,1\right\} }T_{\pmb{X}}\left(n\right)\right)\right)=1+\sum_{n\in\mathbb{N}\backslash\left\{ 0,1\right\} }T_{\pmb{X}}\left(n\right)\end{equation}
so that the Formal Tower Series $\left(1+\sum_{n\in\mathbb{N}\backslash\left\{ 0,1\right\} }T_{\pmb{X}}\left(n\right)\right)$
is invariant under the action of the operator\begin{equation}
\prod_{k\in\mathbb{N}\backslash\left\{ 0\right\} }\left(1+\pmb{\mathfrak{R}}\left(x_{k},\;\centerdot\right)\right)\end{equation}
To get a sense of how such an invariance principle could naturally
arises we consider the function.\begin{equation}
g(x)=1+x+x^{x}+x^{x^{x}}+x^{x^{x^{x}}}+\cdots\end{equation}
and use it to induce a sequence of functions on the vector $\pmb{X}$
who's initial element is\begin{equation}
G_{0}\left(\pmb{X}\right)=\prod_{1\leq k\leq\dim\{\pmb{X}\}}\left(1+g\left(x_{k}\right)\right)\end{equation}
and the other elements of the sequence are defined by the following
recursion\begin{equation}
G_{n+1}\left(\pmb{X}\right)=\left(\prod_{1\leq k\leq\dim\{\pmb{X}\}}\left(1+\pmb{\mathfrak{R}}\left(x_{k},G_{n}\left(\pmb{X}\right)\right)\right)\right)\end{equation}
So that the fundamental invariance principle is re-casted as \begin{equation}
\lim_{n\rightarrow\infty}\left\{ G_{n}\left(\pmb{X}\right)\right\} =1+\sum_{n\in\mathbb{N}\backslash\left\{ 0,1\right\} }T_{\pmb{X}}\left(n\right)\end{equation}
The invariance principle follows from the F. T. A where $\pmb{X}$
will represent the vector whose entries are the distinct primes arranged
in increasing order%
\footnote{the ordering is not necessary for the invariance principle it suffice
to have distinct the entries in $\pmb{X}$ of must corresponding to
distinct primes.%
}. \\
\\
Algorithms such as Buchberger's algorithm in commutative algebra
emphasizes the importance of totally ordering monomials. In our discussion
we shall use the integer ordering to induce a natural ordering on
the towers. Once the towers are totally ordered it becomes rather
straight forward to discuss\textit{ Tower Arithmetic}. Let us encapsulate
the ordering of towers into a metric function $d\left(\centerdot,\centerdot\right)$
so as to embed towers into a metric space $\left(\mathcal{T}\left(\pmb{X}\right),\, d\left(\centerdot,\centerdot\right)\right)$.
The metric space $\left(\mathcal{T}\left(\pmb{X}\right),\, d\left(\centerdot,\centerdot\right)\right)$
allows us express addition of towers through the following relation\[
for\quad T_{\pmb{X}}\left(p\right)\ge T_{\pmb{X}}\left(m\right)\]

\begin{equation}
d\left(T_{\pmb{X}}\left(m\right),\: T_{\pmb{X}}\left(p\right)\right)=d\left(0,\: T_{\pmb{X}}\left(n\right)\right)\Leftrightarrow T_{\pmb{X}}\left(m\right)+T_{\pmb{X}}\left(n\right)=T_{\pmb{X}}\left(p\right).\end{equation}
and we use the convention\begin{equation}
T_{\pmb{X}}\left(0\right)=1\quad and\quad T_{\pmb{X}}\left(0\right)+T_{\pmb{X}}\left(0\right)=x_{1}\end{equation}
Furthermore multiplication of tower follows immediately from the definition
of addition and it is expressed by\[
\left(\prod_{1\le k\le\dim\left\{ \pmb{X}\right\} }x_{k}^{T_{\pmb{X}}\left(m_{k}\right)}\right)\times\left(\prod_{1\le t\le\dim\left\{ \pmb{X}\right\} }x_{t}^{T_{\pmb{X}}\left(n_{t}\right)}\right)=\]
\begin{equation}
\left(\prod_{1\le k\le\dim\left\{ \pmb{X}\right\} }x_{k}^{\left(T_{\pmb{X}}\left(m_{k}\right)+T_{\pmb{X}}\left(n_{k}\right)\right)}\right)\end{equation}
In summary the base level of the product of tower is the union of
the base level of the towers being multiplied while the powers of
corresponding pillars are added.\\
\\
We now consider the case of finite dimensional vectors. Let $\pmb{P}$
be a finite dimensional vector whose entries are made of the smallest
$\dim\left\{ \pmb{P}\right\} $ distinct primes. For convenience we
arrange the primes in increasing order as entries of $\pmb{P}$ we
have \begin{equation}
\pmb{P}=\left(p_{1},\cdots,p_{\dim\left\{ \pmb{P}\right\} }\right)\end{equation}
we define an analogous sequence of functions of $\pmb{P}$ defined
by 

\begin{equation}
G_{0}\left(\pmb{P}\right)=\prod_{1\leq k\leq\dim\left\{ \pmb{P}\right\} }\left(1+g\left(p_{k}\right)\right)\end{equation}
and the recursion\begin{equation}
G_{n+1}\left(\pmb{P}\right)=\prod_{1\leq k\leq\dim\{\pmb{P}\}}\left(1+\pmb{\mathfrak{R}}\left(p_{k},\: G_{n+1}\left(\pmb{P}\right)\right)\right)\end{equation}
in which case we obtain 

\begin{equation}
\lim_{n\rightarrow\infty}\left\{ G_{n}\left(\pmb{P}\right)\right\} =1+\sum_{k\in\mathbb{N}/\{0,1\}}a_{k}\: T_{\pmb{P}}\left(k\right)\end{equation}
where $a_{k}$$\in\left\{ 0,1\right\} $, more specifically $a_{k}=1$
if the tower expansion of $n$ is a tower over $\pmb{P}$ and $a_{k}=0$
otherwise.\\
\\
The preceding discussion raises the following interesting question:
Considering a given finite set of consecutive of integers bounded
by $n$. What is the probability that a number chosen at random contains
a particular prime $p$ in it's tower expansion\\
\\
We now propose a theorem which follows from Euler's argument in
his proof of the Infinity of the primes.\\
\\
\pmb{\emph{Theorem} :} For every finite dimensional vector
$\pmb{P}=\left(p_{1},\cdots,p_{\dim\{\pmb{P}\}}\right)$ whose entries
are made up of distinct primes when considering the limit\begin{equation}
\lim_{n\rightarrow\infty}\left\{ G_{n}\left(\pmb{P}\right)\right\} =1+\sum_{n\in\mathbb{N}/\{0,1\}}a_{n}\; T_{\pmb{P}}\left(n\right).\end{equation}
we have

\begin{equation}
\left(1+\sum_{n\in\mathbb{N}/\{0,1\}}a_{n}\;\left(T_{\pmb{P}}\left(n\right)\right)^{-1}\right)<\infty\end{equation}
\pmb{ Proof }: The convergence follows immediately from the fact
that 

\begin{equation}
\left(1+\sum_{n\in\mathbb{N}/\{0,1\}}a_{n}\;\left(T_{\pmb{P}}\left(n\right)\right)^{-1}\right)<\left(\prod_{1\leq k\leq\dim\{\pmb{P}\}}\left(1-p_{k}{}^{-1}\right)\right)<\infty\end{equation}
\\
\\
The preceding theorem suggests that the rational numbers must
not be too far out of our reach once we are equipped with a concrete
description of the integers as towers. we recall that for a vector
$\pmb{P}=\left(p_{1},\cdots,p_{\dim\{\pmb{P}\}}\right)$ whose entries
are made up of distinct primes. We consider the sequence

\begin{equation}
G_{0}\left(\pmb{P}\right)=\prod_{1\leq k\leq\dim\{\pmb{P}\}}\left(1+g\left(p_{k}\right)\right)\end{equation}

\begin{equation}
G_{n+1}\left(\pmb{P}\right)=\prod_{1\leq k\leq\dim\{\pmb{P}\}}\left(1+\pmb{\mathfrak{R}}\left(p_{k},\: G_{n}\left(\pmb{P}\right)\right)\right)\end{equation}
This sequence may be used to induce the sequence $H_{n}$ defined
by

\begin{equation}
H_{n}\left(\pmb{P}\right)=\prod_{1\leq k\leq\dim\{\pmb{P}\}}\left(\pmb{\mathfrak{R}}\left(p_{k}{}^{-1},\: G_{n}\left(\pmb{P}\right)\right)+1+\text{ }\pmb{\mathfrak{R}}\left(p_{k},G_{n}\left(\pmb{P}\right)\right)\right)\end{equation}
So that the terms in the resulting expression are given by

\begin{equation}
\text{Lim}_{n\rightarrow\infty}\left\{ H_{n}\left(\pmb{P}\right)\right\} =1+\sum_{q\in\mathbb{Q}\backslash\left\{ 0,1\right\} }a_{q}\; T_{\pmb{P}}\left(q\right).\end{equation}
If we seek the complete bijection with the rational we would have
started with the infinite dimensional vector $\pmb{X}$ instead and
considered the following sequence \begin{equation}
G_{0}\left(\pmb{X}\right)=\prod_{1\leq k\leq\dim\{\pmb{X}\}}\left(1+g\left(x_{k}\right)\right)\end{equation}

\begin{equation}
G_{n+1}\left(\pmb{X}\right)=\prod_{1\leq k\leq\dim\{\pmb{X}\}}\left(1+\pmb{\mathfrak{R}}\left(p_{k},\: G_{n}\left(\pmb{X}\right)\right)\right)\end{equation}
This sequence may be used to induce the sequence $H_{n}$ defined
by

\begin{equation}
H_{n}\left(\pmb{X}\right)=\prod_{1\leq k\leq\dim\{\pmb{X}\}}\left(\pmb{\mathfrak{R}}\left(p_{k}{}^{-1},\: G_{n}\left(\pmb{X}\right)\right)+1+\text{ }\pmb{\mathfrak{R}}\left(p_{k},\: G_{n}\left(\pmb{X}\right)\right)\right)\end{equation}
towers in bijections set $\mathbb{Q}\backslash\left\{ 0,1\right\} $
are described by 

\begin{equation}
\text{Lim}_{n\rightarrow\infty}\left\{ H_{n}\left(\pmb{P}\right)\right\} =1+\sum_{q\in\mathbb{Q}\backslash\{0,1\}}T_{\pmb{P}}\left(q\right).\end{equation}

\subsection{Tower Sieve Algorithm }

Sieves play an important role in Number theory, We propose to investigate
here a novel sieve algorithm based on the arithmetic of towers. Let
$\pmb{P}=\left(p_{1},\cdots,p_{\dim\{\pmb{P}\}}\right)$ denote vector
whose entries are all the primes less the $2^{t}$ for some $1\leq t$,
the goal is to determine the primes in the range $\left[2^{t},2^{1+t}\right]$.
The algorithm consists in computing the following recursion

\begin{equation}
\left.\left.g_{s}(x)=1+x+x^{x}+\text{ }\cdots\text{ }+\left(x^{x^{.\cdot{}^{\cdot}{}^{x}}}\right.\right\} \text{last }\text{term }\text{of }\text{height }s\right)\end{equation}

\begin{equation}
G_{0}\left(\pmb{P}\right)=\prod_{1\leq k\leq\dim\{\pmb{P}\}}g_{s_{k}}\left(p_{k}\right)\end{equation}
where $p_{k}{}^{p_{k}{}^{.\cdot{}^{\cdot}{}^{p_{k}}}}$is the last
term in the expression $g_{s_{k}}\left(p_{k}\right)$ is such that
$p_{k}{}^{p_{k}{}^{.\cdot{}^{\cdot}{}^{p_{k}}}}\ll2^{1+t}$ \begin{equation}
G_{n+1}\left(\pmb{P}\right)=\left(\prod_{1\leq k\leq\dim\{\pmb{P}\}}\left(1+\pmb{\mathfrak{R}}\left(p_{k},\: G_{n}\left(\pmb{P}\right)\right)\right)\right)\end{equation}
we stop the recursion at the $m^{\text{th}}$ iteration if all the
towers remaining in the polytower difference $\left(G_{m+1}\left(\pmb{P}\right)-G_{m}\left(\pmb{P}\right)\right)$
are are towers greater than $2^{1+t}$. \\
\\
At the heart of the recursive algorithm is the fact that the recursion
determines the tower expansions of integers in the interval $\left[2^{t},\:2^{1+t}\right]$
with the exception of towers expansion which contains primes which
are not less than $2^{t}$. Furthermore assuming that we order the
towers the gaps of size $2$ in the list determines the exact location
of primes in the range of interest. This provide a constructive proof
of the fact that there is always at least one prime in the range $\left[2^{t},\,2^{1+t}\right]$
for any value of $1\leq t$. The sieves algorithm we discussed above
is rather different from Eratosthenes sieve in that it generates the
composite and indicates exactly where the primes ought to be found
and more importantly as oppose to some variants of Eratosthenes sieve
methods our algorithm ensures that each composite is generated exactly
once. Let us briefly go through the steps the algorithm with the case
$t=2$ for illustration purposes

\subsection{Illustration of the Algorithm}

We illustrate the algorithm using Mathematica. \\
Let $\pmb{X}=\left(x_{1},x_{2}\right)$ the Mathematica commands
used are \[
\pmb{g_{8}\text{:=}1+x_{k}}\]
\[
\pmb{\text{For}\left[i=1,i<7,i\text{++},g_{8}=\left(1+\text{Total}\left[\left(x_{k}{}^{\wedge}\text{List}\text{@@}(g)\right)\right]\right)\right]}\]
from which $\pmb{g_{8}}$ is given by 

\noindent \begin{equation}
\pmb{g}=1+x_{k}+x_{k}^{x_{k}}+x_{k}^{x_{k}^{x_{k}}}+x_{k}^{x_{k}^{x_{k}^{x_{k}}}}+x_{k}^{x_{k}^{x_{k}^{x_{k}^{x_{k}}}}}+x_{k}^{x_{k}^{x_{k}^{x_{k}^{x_{k}^{x_{k}}}}}}+x_{k}^{x_{k}^{x_{k}^{x_{k}^{x_{k}^{x_{k}^{x_{k}}}}}}}\end{equation}
The Mathematica commands for the recursion are given by \[
\pmb{G_{0}\text{:=}1}\]
\[
\pmb{\text{For}\left[k=1,k<3,k\text{++},G_{0}=\text{Expand}\left[G_{0}\left(1+x_{k}\right)\right]\right]}\]
after these first commands we have \begin{equation}
\pmb{G_{0}}=1+x_{1}+x_{2}+x_{1}x_{2}\end{equation}
Here is an overview of the typical intermediary steps required to
compute the recursion.

\noindent \[
\pmb{\text{List}\text{@@}\left(G_{0}\right)}=\left\{ 1,x_{1},x_{2},x_{1}x_{2}\right\} \]

\noindent \[
\pmb{\left\{ 1,x_{1},x_{2},x_{1}x_{2}\right\} \text{/.}x_{1}\to2\text{/.}x_{2}\to3}=\left\{ 1,2,3,6\right\} \]
\[
\pmb{\left(x_{1}{}^{\wedge}\text{List}\text{@@}\left(G_{0}\right)\right)}=\left\{ x_{1},x_{1}^{x_{1}},x_{1}^{x_{2}},x_{1}^{x_{1}x_{2}}\right\} \]
\[
\pmb{\text{Total}\left[\left(x_{1}{}^{\wedge}\text{List}\text{@@}\left(G_{0}\right)\right)\right]}=x_{1}+x_{1}^{x_{1}}+x_{1}^{x_{2}}+x_{1}^{x_{1}x_{2}}\]
\[
\pmb{\text{Expand}\left[\left(1+\text{Total}\left[\left(x_{1}{}^{\wedge}\text{List}\text{@@}\left(G_{0}\right)\right)\right]\right)\left(1+\text{Total}\left[\left(x_{2}{}^{\wedge}\text{List}\text{@@}\left(1+x_{1}\right)\right)\right]\right)\right]}\]

\noindent \[
=1+x_{1}+x_{1}^{x_{1}}+x_{1}^{x_{2}}+x_{1}^{x_{1}x_{2}}+x_{2}+x_{1}x_{2}+x_{1}^{x_{1}}x_{2}+x_{1}^{x_{2}}x_{2}+\]
\[
x_{1}^{x_{1}x_{2}}x_{2}+x_{2}^{x_{1}}+x_{1}x_{2}^{x_{1}}+x_{1}^{x_{1}}x_{2}^{x_{1}}+x_{1}^{x_{2}}x_{2}^{x_{1}}+x_{1}^{x_{1}x_{2}}x_{2}^{x_{1}}\]
The polytower contains the towers of interest and the list bellow
depicts the corresponding numbers.\[
\pmb{\text{L0}=\text{Sort}[\left[\text{List}\text{@@}\text{Expand}\left[\left(1+\text{Total}\left[\left(x_{1}{}^{\wedge}\text{List}\text{@@}\left(\text{Expand}\left[\left(1+x_{1}\right)\left(1+x_{2}\right)\right]\right)\right)\right]\right)\right.\right.}\]
\[
\pmb{\left.\left.\left(1+\text{Total}\left[\left(x_{2}{}^{\wedge}\text{List}\text{@@}\left(1+x_{1}\right)\right)\right]\right)\right]/.x_{1}\to2\text{/.}x_{2}\to3\right]}\]
\begin{equation}
=\{1,2,3,4,6,8,9,12,18,24,36,64,72,192,576\}\end{equation}
If we add in the primes determined by the list we get the following
sequence of towers and their corresponding numbers.\[
\pmb{F\text{:=}}\]
\[
\pmb{\text{Expand}\left[\left(1+\text{Total}\left[\left(x_{1}{}^{\wedge}\text{List}\text{@@}\left(\text{Expand}\left[\left(1+x_{1}\right)\left(1+x_{2}\right)\right]\right)\right)\right]\right)\left(1+\text{Total}\left[\left(x_{2}{}^{\wedge}\text{List}\text{@@}\left(1+x_{1}\right)\right)\right]\right)\right.}\]
\[
\pmb{\left.\left(1+x_{3}\right)\left(1+x_{4}\right)\right]}\]
if we list the resulting towers we get \[
\pmb{\text{List}\text{@@}F}=\]

\noindent \[
\left\{ 1,x_{1},x_{1}^{x_{1}},x_{1}^{x_{2}},x_{1}^{x_{1}x_{2}},x_{2},x_{1}x_{2},x_{1}^{x_{1}}x_{2},x_{1}^{x_{2}}x_{2},x_{1}^{x_{1}x_{2}}x_{2},x_{2}^{x_{1}},x_{1}x_{2}^{x_{1}},\right.\]

\noindent \[
x_{1}^{x_{1}}x_{2}^{x_{1}},x_{1}^{x_{2}}x_{2}^{x_{1}},x_{1}^{x_{1}x_{2}}x_{2}^{x_{1}},x_{3},x_{1}x_{3},x_{1}^{x_{1}}x_{3},x_{1}^{x_{2}}x_{3},x_{1}^{x_{1}x_{2}}x_{3},x_{2}x_{3},x_{1}x_{2}x_{3},\]
\[
x_{1}^{x_{1}}x_{2}x_{3},x_{1}^{x_{2}}x_{2}x_{3},x_{1}^{x_{1}x_{2}}x_{2}x_{3},x_{2}^{x_{1}}x_{3},x_{1}x_{2}^{x_{1}}x_{3},x_{1}^{x_{1}}x_{2}^{x_{1}}x_{3},x_{1}^{x_{2}}x_{2}^{x_{1}}x_{3},\]
\[
x_{1}^{x_{1}x_{2}}x_{2}^{x_{1}}x_{3},x_{4},x_{1}x_{4},x_{1}^{x_{1}}x_{4},x_{1}^{x_{2}}x_{4},x_{1}^{x_{1}x_{2}}x_{4},x_{2}x_{4},x_{1}x_{2}x_{4},x_{1}^{x_{1}}x_{2}x_{4},\]
\[
x_{1}^{x_{2}}x_{2}x_{4},x_{1}^{x_{1}x_{2}}x_{2}x_{4},x_{2}^{x_{1}}x_{4},x_{1}x_{2}^{x_{1}}x_{4},x_{1}^{x_{1}}x_{2}^{x_{1}}x_{4},x_{1}^{x_{2}}x_{2}^{x_{1}}x_{4},x_{1}^{x_{1}x_{2}}x_{2}^{x_{1}}x_{4}\]
\[
,x_{3}x_{4},x_{1}x_{3}x_{4},x_{1}^{x_{1}}x_{3}x_{4},x_{1}^{x_{2}}x_{3}x_{4},x_{1}^{x_{1}x_{2}}x_{3}x_{4},x_{2}x_{3}x_{4},x_{1}x_{2}x_{3}x_{4},\]
\[
x_{1}^{x_{1}}x_{2}x_{3}x_{4},x_{1}^{x_{2}}x_{2}x_{3}x_{4},x_{1}^{x_{1}x_{2}}x_{2}x_{3}x_{4},x_{2}^{x_{1}}x_{3}x_{4},\]

\noindent \begin{equation}
\left.x_{1}x_{2}^{x_{1}}x_{3}x_{4},x_{1}^{x_{1}}x_{2}^{x_{1}}x_{3}x_{4},x_{1}^{x_{2}}x_{2}^{x_{1}}x_{3}x_{4},x_{1}^{x_{1}x_{2}}x_{2}^{x_{1}}x_{3}x_{4}\right\} \end{equation}
The corresponding list of integer is determined by the following Mathematica
commands

\[
\pmb{\text{L1}=\text{Sort}\left[\text{List}\text{@@}F\text{/.}x_{1}\to2\text{/.}x_{2}\to3\text{/.}x_{3}\to5\text{/.}x_{4}\to7\right]}\]
which yields the following list of integers.

\[
\left\{ 1,2,3,4,5,6,7,8,9,10,12,14,15,18,20,21,24,28,30,35,\right.\]
\[
36,40,42,45,56,60,63,64,70,72,84,90,105,120,126,140,168,180,\]
\[
192,210,252,280,315,320,360,420,448,504,576,630,840\]

\begin{equation}
\left.,960,1260,1344,2240,2520,2880,4032,6720,20160\right\} \end{equation}

We note that the simple insertion of these towers determines the prime
number 11 and 13. The next step is to discuss a slight modification
of the straightforward technique discussed above . We see from the
algorithm that we are very strongly incentives to reduce the number
of terms appearing in the expression for instance it is clear that
the number should not be included so we subs-tract it and consider
the following expression . 

\noindent \[
\pmb{A\text{:=}1+\text{Total}\left[\left(x_{1}{}^{\wedge}\text{List}\text{@@}\left(\text{Expand}\left[\left(1+x_{1}\right)\left(1+x_{2}\right)\right]-x_{1}x_{2}\right)\right)\right]}\]
The previous would correspond to what I would refer to as a renormalization
step\[
\pmb{B\text{:=}1+x_{2}}\]
The list of towers and the corresponding list of numbers generated
by the reduced product is given by performing by re-normalizing as
follows 

\noindent \[
\pmb{\text{Expand}[AB]-x_{1}^{x_{1}}x_{2}-x_{1}^{x_{2}}x_{2}}\]
\begin{equation}
=1+x_{1}+x_{1}^{x_{1}}+x_{1}^{x_{2}}+x_{2}+x_{1}x_{2}\end{equation}

\noindent \[
\pmb{\text{Sort}\left[\text{List}\text{@@}\left(\text{Expand}[AB]-x_{1}^{x_{1}}x_{2}-x_{1}^{x_{2}}x_{2}\right)\text{/.}x_{1}\to2\text{/.}x_{2}\to3\right]}\]

\noindent \begin{equation}
=\{1,2,3,4,6,8\}\end{equation}
Now if we want to improve the estimate of the sum of the primes in
the range and than we must remove from the sum the towers associated
with the number in the range of 1 and as follows 

\[
\pmb{\text{Expand}\left[\left(1+\text{Total}\left[\left(x_{1}{}^{\wedge}\text{List}\text{@@}\left(\text{Expand}\left[\left(1+x_{1}\right)\left(1+x_{2}\right)\right]\right)\right)\right]-x_{1}^{x_{1}x_{2}}\right)\left(1+x_{2}\right)\right]-x_{1}^{x_{1}}x_{2}-}\]
\[
\pmb{x_{1}^{x_{2}}x_{2}-1-x_{1}-x_{2}}\]

\noindent \begin{equation}
=x_{1}^{x_{1}}+x_{1}^{x_{2}}+x_{1}x_{2}\end{equation}
Here is now the estimate if the sum of the primes in the range , is
given by the following 

\noindent \[
\pmb{G\text{:=}\text{Expand}\left[AB\left(1+x_{3}\right)\left(1+x_{4}\right)\right]}\]

\noindent \[
\pmb{\text{L2}=\text{Sort}\left[\text{List}\text{@@}G\text{/.}x_{1}\to2\text{/.}x_{2}\to3\text{/.}x_{3}\to5\text{/.}x_{4}\to7\right]}\]
\[
=\left\{ 1,2,3,4,5,6,7,8,10,12,14,15,20,21,24,28,30,35,40,\right.\]
\begin{equation}
\left.42,56,60,70,84,105,120,140,168,210,280,420,840\right\} \end{equation}
We note that the impact of the reduction is significant we have almost
reduced the length we are considering by a factor of 2 as depicted
bellow which in of itself is pretty amazing.

\noindent \[
\pmb{\text{Dimensions}[\text{L1}]}\]
\begin{equation}
=\{60\}\end{equation}
\[
\pmb{\text{Dimensions}[\text{L2}]}\]

\noindent \begin{equation}
\{32\}\end{equation}
if one just wants to generate the trees without any concern for the
the renormalization business here is how one proceeds.\[
\pmb{G_{0}\text{:=}1}\]

\[
\pmb{\text{For}\left[k=1,k<3,k\text{++},G_{0}=\text{Expand}\left[G_{0}\left(1+x_{k}\right)\right]\right]}\]
\[
\pmb{G_{1}\text{:=}1}\pmb{\text{For}\left[k=1,k<3,k\text{++},G_{1}=\text{Expand}\left[G_{1}\left(1+\text{Total}\left[\left(x_{k}{}^{\wedge}\text{List}\text{@@}\left(G_{0}\right)\right)\right]\right)\right]\right]}\pmb{G_{1}}\]

\noindent \[
=\left\{ 1+x_{1}+x_{1}^{x_{1}}+x_{1}^{x_{2}}+x_{1}^{x_{1}x_{2}}+x_{2}+x_{1}x_{2}+x_{1}^{x_{1}}x_{2}+x_{1}^{x_{2}}x_{2}+\right.\]
\[
x_{1}^{x_{1}x_{2}}x_{2}+x_{2}^{x_{1}}+x_{1}x_{2}^{x_{1}}+x_{1}^{x_{1}}x_{2}^{x_{1}}+x_{1}^{x_{2}}x_{2}^{x_{1}}+x_{1}^{x_{1}x_{2}}x_{2}^{x_{1}}+x_{2}^{x_{2}}+x_{1}x_{2}^{x_{2}}+\]
\begin{equation}
\left.x_{1}^{x_{1}}x_{2}^{x_{2}}+x_{1}^{x_{2}}x_{2}^{x_{2}}+x_{1}^{x_{1}x_{2}}x_{2}^{x_{2}}+x_{2}^{x_{1}x_{2}}+x_{1}x_{2}^{x_{1}x_{2}}+x_{1}^{x_{1}}x_{2}^{x_{1}x_{2}}+x_{1}^{x_{2}}x_{2}^{x_{1}x_{2}}+x_{1}^{x_{1}x_{2}}x_{2}^{x_{1}x_{2}}\right\} \end{equation}

\noindent which results in the following list of integers.\[
\{1,2,3,4,6,8,9,12,18,24,27,36,54,64,72,108,192\]

\noindent \begin{equation}
,216,576,729,1458,1728,2916,5832,46656\}\end{equation}
The preceding sequence of numbers allowed us to determine that the
primes 5 and 7 are missing from the list.

\subsection{Recovering the Ordering.}

What we would presumably like would be to recover the ordering of
the integers by simply manipulating the tower expansion. This would
also mean that we would recover the values of the primes. this would
not be an easy task but it has the merit of showing how in a sense
all the primes are determined by the first two integers. We therefore
define

\begin{equation}
g_{s_{i}}\left(x_{i,t}\right)=1+x_{i,t}+x_{i,t}{}^{x_{i,t}}+\text{ }\cdots\text{ }+\left(\left.x_{i,t}{}^{x_{i,t}{}^{.\cdot{}^{\cdot}{}^{x_{i,t}}}}\right\} \text{last }\text{term }\text{of }\text{height }s_{i}\right)\end{equation}

\begin{equation}
G_{0}\left(\pmb{X}\right)=\prod_{1\leq i\leq\dim\{\pmb{X}\}}\left(1+g_{s_{i}}\left(x_{i,1}\right)\right)\end{equation}
We then define the following recurrence

\begin{equation}
G_{m}\left(\pmb{X}\right)=\prod_{1\leq i\leq\dim\{\pmb{X}\}}\left(1+\pmb{\mathfrak{R}}\left(x_{i,m},\; G_{m-1}\left(\pmb{X}\right)\right)\right)\end{equation}
The basic idea is to perform substitution of the variables in the
towers so as to end up with polynomials in the single variable $x$
which can in turn be ordered by recovering the binary decimal expansion
corresponding to the numbers. Note that the polynomials obtained when
evaluated at $x=2$ give the value of the integers that the number
represents. Let us illustrate the computation for $\pmb{X}=\left(x_{1},x_{2}\right)$.
Here are the Mathematica Commands illustrating the overall approach

\noindent \[
\pmb{G_{0}=\text{Expand}\left[\left(1+x_{1,1}\right)\left(1+x_{2,1}\right)\right]}\]

\noindent \begin{equation}
=1+x_{1,1}+x_{2,1}+x_{1,1}x_{2,1}\end{equation}
\[
\pmb{\text{Expand}\left[\left(1+x_{1,1}\right)\left(1+x_{2,1}\right)\right]}\]
\begin{equation}
1+x_{1,1}+x_{2,1}+x_{1,1}x_{2,1}\end{equation}
\[
\pmb{G_{1}=}\]
\[
\pmb{\text{Expand}\left[\left(1+\text{Total}\left[\left(x_{1,2}{}^{\wedge}\text{List}\text{@@}\left(\text{Expand}\left[\left(1+x_{1,1}\right)\left(1+x_{2,1}\right)\right]-x_{1,1}x_{2,1}-x_{2,1}\right)\right)\right]\right)\right.}\]
\[
\pmb{\left.\left(1+x_{2,1}\right)\right]-x_{1,2}^{x_{1,1}}x_{2,1}}\]
\begin{equation}
=1+x_{1,2}+x_{1,2}^{x_{1,1}}+x_{2,1}+x_{1,2}x_{2,1}\end{equation}
note that we have performed the renormalization in the expressions
above to obtain the following polynomial expressions.

\noindent \begin{equation}
\pmb{\text{Expand}\left[\text{List}\text{@@}G_{0}\text{/.}x_{1,1}\to2\text{/.}x_{1,2}\to x\text{/.}x_{2,1}\to1+x\right]}\end{equation}

\noindent \begin{equation}
=\left\{ 1,x,x^{2},1+x,x+x^{2}\right\} \end{equation}
The expression $G_{1}$ determines all the primes between $x^{2}$
and $x^{2}+x$ so as to determine the list of consecutive integers
from $1$ to $x^{2}+x$ it suffice to add $x_{3,1}$ to the $G_{1}$
as follows

\noindent \[
\pmb{T_{1}=G_{1}+x_{3,1}}\]

\noindent \begin{equation}
=1+x_{1,2}+x_{1,2}^{x_{1,1}}+x_{2,1}+x_{1,2}x_{2,1}+x_{3,1}\end{equation}

\noindent \[
\pmb{\text{Expand}\left[\text{List}\text{@@}T_{1}\text{/.}x_{1,1}\to2\text{/.}x_{1,2}\to x\text{/.}x_{2,1}\to1+x\text{/.}x_{3,1}\to x^{2}+1\right]}\]

\noindent \begin{equation}
=\left\{ 1,x,x^{2},1+x,x+x^{2},1+x^{2}\right\} \end{equation}
\[
\pmb{R_{2}=\text{Expand}\left[\left(1+\text{Total}\left[\left(x_{1,2}{}^{\wedge}\text{List}\text{@@}\left(\text{Expand}\left[\left(1+x_{1,1}\right)\left(1+x_{2,1}\right)\right]-x_{1,1}x_{2,1}\right)\right)\right]\right)\right.}\]

\noindent \[
\pmb{\left.\left(1+\text{Total}\left[\left(x_{2,2}{}^{\wedge}\text{List}\text{@@}\left(\text{Expand}\left[\left(1+x_{1,1}\right)\left(1+x_{2,1}\right)\right]-x_{1,1}x_{2,1}-x_{2,1}\right)\right)\right]\right)\left(1+x_{3,1}\right)\right]}\]
\[
1+x_{1,2}+x_{1,2}^{x_{1,1}}+x_{1,2}^{x_{2,1}}+x_{2,2}+x_{1,2}x_{2,2}+x_{1,2}^{x_{1,1}}x_{2,2}+x_{1,2}^{x_{2,1}}x_{2,2}+x_{2,2}^{x_{1,1}}+\]
\[
x_{1,2}x_{2,2}^{x_{1,1}}+x_{1,2}^{x_{1,1}}x_{2,2}^{x_{1,1}}+x_{1,2}^{x_{2,1}}x_{2,2}^{x_{1,1}}+x_{3,1}+x_{1,2}x_{3,1}+x_{1,2}^{x_{1,1}}x_{3,1}+\]
\[
x_{1,2}^{x_{2,1}}x_{3,1}+x_{2,2}x_{3,1}+x_{1,2}x_{2,2}x_{3,1}+x_{1,2}^{x_{1,1}}x_{2,2}x_{3,1}+x_{1,2}^{x_{2,1}}x_{2,2}x_{3,1}+\]

\noindent \begin{equation}
x_{2,2}^{x_{1,1}}x_{3,1}+x_{1,2}x_{2,2}^{x_{1,1}}x_{3,1}+x_{1,2}^{x_{1,1}}x_{2,2}^{x_{1,1}}x_{3,1}+x_{1,2}^{x_{2,1}}x_{2,2}^{x_{1,1}}x_{3,1}\end{equation}

\[
\pmb{G_{2}=}\pmb{R_{2}-\left(x_{1,2}^{x_{1,1}}x_{3,1}+x_{1,2}^{x_{2,1}}x_{3,1}+x_{2,2}x_{3,1}+x_{1,2}x_{2,2}x_{3,1}+\right.}\]
\[
\pmb{x_{1,2}^{x_{1,1}}x_{2,2}x_{3,1}+x_{1,2}^{x_{2,1}}x_{2,2}x_{3,1}+x_{1,2}^{x_{2,1}}x_{2,2}+x_{1,2}x_{2,2}^{x_{1,1}}x_{3,1}+}\]

\noindent \begin{equation}
\pmb{\left.x_{1,2}^{x_{1,1}}x_{2,2}^{x_{1,1}}x_{3,1}+x_{1,2}^{x_{2,1}}x_{2,2}^{x_{1,1}}x_{3,1}+x_{2,2}^{x_{1,1}}x_{3,1}+x_{1,2}^{x_{1,1}}x_{2,2}^{x_{1,1}}+x_{1,2}^{x_{2,1}}x_{2,2}^{x_{1,1}}+x_{1,2}x_{2,2}^{x_{1,1}}\right)}\end{equation}

\noindent \begin{equation}
=1+x_{1,2}+x_{1,2}^{x_{1,1}}+x_{1,2}^{x_{2,1}}+x_{2,2}+x_{1,2}x_{2,2}+x_{1,2}^{x_{1,1}}x_{2,2}+x_{2,2}^{x_{1,1}}+x_{3,1}+x_{1,2}x_{3,1}\end{equation}

\noindent \[
\pmb{\text{Expand}\left[\text{List}\text{@@}G_{2}\text{/.}x_{1,2}\to x\text{/.}x_{1,1}\to2\right.}\]
\[
\pmb{\left.\text{/.}x_{2,1}\to3\text{/.}x_{2,2}\to(1+x)\text{/.}x_{3,1}\to\left(x^{2}+1\right)\right]}\]

\noindent \begin{equation}
=\left\{ 1,x,x^{2},x^{3},1+x,x+x^{2},x^{2}+x^{3},1+2x+x^{2},1+x^{2},x+x^{3}\right\} \end{equation}

\noindent \[
\pmb{\text{Sort}\left[\text{Expand}\left[\text{List}\text{@@}G_{2}\text{/.}x_{1,2}\to x\text{/.}x_{1,1}\to2\text{/.}x_{2,1}\to3\right.\right.}\]
\[
\pmb{\left.\left.\text{/.}x_{2,2}\to(1+x)\text{/.}x_{3,1}\to\left(x^{2}+1\right)\right]\text{/.}x\to2\right]}\]

\noindent \begin{equation}
=\{1,2,3,4,5,6,8,9,10,12\}\end{equation}

\noindent \[
\pmb{T_{2}=G_{2}+x_{4,1}+x_{5,1}}=\]

\noindent \begin{equation}
1+x_{1,2}+x_{1,2}^{x_{1,1}}+x_{1,2}^{x_{2,1}}+x_{2,2}+x_{1,2}x_{2,2}+x_{1,2}^{x_{1,1}}x_{2,2}+x_{2,2}^{x_{1,1}}+x_{3,1}+x_{1,2}x_{3,1}+x_{4,1}+x_{5,1}\end{equation}
In the following list we express the tower expansion for the consecutive
integers in the range 0 and $p_{5}$.\[
\pmb{\text{Expand}\left[\text{List}\text{@@}T_{1}\text{/.}x_{1,2}\to p_{1}\text{/.}x_{1,1}\to p_{1}\right.}\]
\[
\pmb{\left.\text{/.}x_{2,1}\to p_{2}\text{/.}x_{2,2}\to p_{2}\text{/.}x_{3,1}\to p_{3}\text{/.}x_{4,1}\to p_{4}\text{/.}x_{5,1}\to p_{5}\right]}\]

\noindent \begin{equation}
=\left\{ 1,p_{1},p_{1}^{p_{1}},p_{1}^{p_{2}},p_{2},p_{1}p_{2},p_{1}^{p_{1}}p_{2},p_{2}^{p_{1}},p_{3},p_{1}p_{3},p_{4},p_{5}\right\} \end{equation}
In the following list we express the tower expansion for the consecutive
integers in the range 0 and $p_{5}$ in the special polynomial form.\[
\pmb{\text{Expand}\left[\text{List}\text{@@}T_{1}\text{/.}x_{1,2}\to x\text{/.}x_{1,1}\to2\text{/.}x_{2,1}\to3\text{/.}x_{2,2}\to(1+x)\right.}\]

\noindent \[
\pmb{\left.\text{/.}x_{3,1}\to\left(x^{2}+1\right)\text{/.}x_{4,1}\to x^{2}+x+1\text{/.}x_{5,1}\to x^{3}+x^{2}+1\right]}\]

\noindent \begin{equation}
=\left\{ 1,x,x^{2},x^{3},1+x,x+x^{2},x^{2}+x^{3},1+2x+x^{2},1+x^{2},x+x^{3},1+x+x^{2},1+x^{2}+x^{3}\right\} \end{equation}
One continues the process illustrated above to obtain The polynomial
representation of the primes in the sought after range. Furthermore
the above expansion suggest a alternative way of representing numbers
in the sense that these are reducible polynomials in the field $\left[\left\{ +,\times\right\} ;\left\{ 0,1\right\} \right]$
and the factorization can be recovered in polynomial time. The problems
of course is that addition completely messes things up. So in order
to recover correct polynomial representation of a given integers after
performing one or several addition consists in generating towers in
the neighborhood of the sought after polynomials and check that they
both evaluate to the same number. Another important fact following
from the discussion above is the fact that the arithmetic is being
performed in a manner which mimics operation on sets one can think
of these operations as operations on forests with rooted trees with
colored vertices where each one one of the primes represent a color.

\subsection{Getting hold of the rationals.}

We recall that 

\begin{equation}
g\left(x\right)=1+x+x^{x}+x^{x^{x}}+x^{x^{x^{x}}}+x^{x^{x^{x^{x}}}}+\cdots\end{equation}
Let \begin{equation}
\boldsymbol{P}\equiv\left(p_{1},\cdots,p_{k},\cdots\right)\end{equation}
denote the vectors of the primes.

\begin{equation}
\lim_{n\rightarrow\infty}\left\{ H_{n}\left(\boldsymbol{P}\right)=\prod_{1\le k\le|\boldsymbol{P}|}\left(\pmb{\mathfrak{R}}\left(p_{k}^{-1},G_{n}\left(\boldsymbol{P}\right)\right)+1+\pmb{\mathfrak{R}}\left(p_{k},G_{n}\left(\boldsymbol{P}\right)\right)\right)\right\} \end{equation}
The terms in the sums are in bijective correspondence with the element
on $\mathbb{Q}$the code for producing rational is given bellow. In
what follows we illustrate how to generate towers associated with
rational numbers.

\[
\pmb{g_{t}\text{:=}1+x_{t}}\]
\[
\pmb{\text{For}\left[i=1,i<1,i\text{++},g_{t}=\left(1+\text{Total}\left[\left(x_{t}{}^{\wedge}\text{List}\text{@@}\left(g_{t}\right)\right)\right]\right)\right]}\]
\[
\pmb{G_{0}\text{:=}1}\]

\[
\pmb{\text{For}\left[k=1,k<3,k\text{++},G_{0}=\text{Expand}\left[G_{0}\left(g_{t}\text{/.}t\to k\right)\right]\right]}\]
\[
\pmb{G_{1}\text{:=}1}\]
\[
\pmb{\text{For}\left[k=1,k<3,k\text{++},G_{1}=\text{Expand}\left[G_{1}\left(1+\text{Total}\left[\left(x_{k}{}^{\wedge}\text{List}\text{@@}\left(G_{0}\right)\right)\right]\right)\right]\right]}\]
\[
\pmb{H_{1}\text{:=}1}\]
\[
\pmb{\text{For}\left[k=1,k<3,k\text{++},\right.}\]
\[
\pmb{\left.H_{1}=\text{Expand}\left[H_{1}\left(\text{Total}\left[\left(x_{k}{}^{-1}\right){}^{\wedge}\text{List}\text{@@}\left(G_{0}\right)\right]+1+\text{Total}\left[x_{k}{}^{\wedge}\text{List}\text{@@}\left(G_{0}\right)\right]\right)\right]\right]}\]
\[
\pmb{H_{1}}\]
\[
=1+\left(\frac{1}{x_{1}}\right){}^{x_{1}}+\left(\frac{1}{x_{1}}\right){}^{x_{2}}+\left(\frac{1}{x_{1}}\right){}^{x_{1}x_{2}}+\frac{1}{x_{1}}+x_{1}+x_{1}^{x_{1}}+x_{1}^{x_{2}}+x_{1}^{x_{1}x_{2}}+\]
\[
\left(\frac{1}{x_{2}}\right){}^{x_{1}}+\left(\frac{1}{x_{1}}\right){}^{x_{1}}\left(\frac{1}{x_{2}}\right){}^{x_{1}}+\left(\frac{1}{x_{1}}\right){}^{x_{2}}\left(\frac{1}{x_{2}}\right){}^{x_{1}}+\left(\frac{1}{x_{1}}\right){}^{x_{1}x_{2}}\left(\frac{1}{x_{2}}\right){}^{x_{1}}+\frac{\left(\frac{1}{x_{2}}\right){}^{x_{1}}}{x_{1}}+\]
\[
x_{1}\left(\frac{1}{x_{2}}\right){}^{x_{1}}+x_{1}^{x_{1}}\left(\frac{1}{x_{2}}\right){}^{x_{1}}+x_{1}^{x_{2}}\left(\frac{1}{x_{2}}\right){}^{x_{1}}+x_{1}^{x_{1}x_{2}}\left(\frac{1}{x_{2}}\right){}^{x_{1}}+\left(\frac{1}{x_{2}}\right){}^{x_{2}}+\left(\frac{1}{x_{1}}\right){}^{x_{1}}\left(\frac{1}{x_{2}}\right){}^{x_{2}}+\]
\[
\left(\frac{1}{x_{1}}\right){}^{x_{2}}\left(\frac{1}{x_{2}}\right){}^{x_{2}}+\left(\frac{1}{x_{1}}\right){}^{x_{1}x_{2}}\left(\frac{1}{x_{2}}\right){}^{x_{2}}+\frac{\left(\frac{1}{x_{2}}\right){}^{x_{2}}}{x_{1}}+x_{1}\left(\frac{1}{x_{2}}\right){}^{x_{2}}+x_{1}^{x_{1}}\left(\frac{1}{x_{2}}\right){}^{x_{2}}+x_{1}^{x_{2}}\left(\frac{1}{x_{2}}\right){}^{x_{2}}+\]
\[
x_{1}^{x_{1}x_{2}}\left(\frac{1}{x_{2}}\right){}^{x_{2}}+\left(\frac{1}{x_{2}}\right){}^{x_{1}x_{2}}+\left(\frac{1}{x_{1}}\right){}^{x_{1}}\left(\frac{1}{x_{2}}\right){}^{x_{1}x_{2}}+\left(\frac{1}{x_{1}}\right){}^{x_{2}}\left(\frac{1}{x_{2}}\right){}^{x_{1}x_{2}}+\left(\frac{1}{x_{1}}\right){}^{x_{1}x_{2}}\left(\frac{1}{x_{2}}\right){}^{x_{1}x_{2}}+\]
\[
\frac{\left(\frac{1}{x_{2}}\right){}^{x_{1}x_{2}}}{x_{1}}+x_{1}\left(\frac{1}{x_{2}}\right){}^{x_{1}x_{2}}+x_{1}^{x_{1}}\left(\frac{1}{x_{2}}\right){}^{x_{1}x_{2}}+x_{1}^{x_{2}}\left(\frac{1}{x_{2}}\right){}^{x_{1}x_{2}}+x_{1}^{x_{1}x_{2}}\left(\frac{1}{x_{2}}\right){}^{x_{1}x_{2}}+\frac{1}{x_{2}}+\frac{\left(\frac{1}{x_{1}}\right){}^{x_{1}}}{x_{2}}+\]
\[
\frac{\left(\frac{1}{x_{1}}\right){}^{x_{2}}}{x_{2}}+\frac{\left(\frac{1}{x_{1}}\right){}^{x_{1}x_{2}}}{x_{2}}+\frac{1}{x_{1}x_{2}}+\frac{x_{1}}{x_{2}}+\frac{x_{1}^{x_{1}}}{x_{2}}+\frac{x_{1}^{x_{2}}}{x_{2}}+\frac{x_{1}^{x_{1}x_{2}}}{x_{2}}+x_{2}+\left(\frac{1}{x_{1}}\right){}^{x_{1}}x_{2}+\left(\frac{1}{x_{1}}\right){}^{x_{2}}x_{2}+\]
\[
\left(\frac{1}{x_{1}}\right){}^{x_{1}x_{2}}x_{2}+\frac{x_{2}}{x_{1}}+x_{1}x_{2}+x_{1}^{x_{1}}x_{2}+x_{1}^{x_{2}}x_{2}+x_{1}^{x_{1}x_{2}}x_{2}+x_{2}^{x_{1}}+\left(\frac{1}{x_{1}}\right){}^{x_{1}}x_{2}^{x_{1}}+\left(\frac{1}{x_{1}}\right){}^{x_{2}}x_{2}^{x_{1}}+\]
\[
\left(\frac{1}{x_{1}}\right){}^{x_{1}x_{2}}x_{2}^{x_{1}}+\frac{x_{2}^{x_{1}}}{x_{1}}+x_{1}x_{2}^{x_{1}}+x_{1}^{x_{1}}x_{2}^{x_{1}}+x_{1}^{x_{2}}x_{2}^{x_{1}}+x_{1}^{x_{1}x_{2}}x_{2}^{x_{1}}+x_{2}^{x_{2}}+\left(\frac{1}{x_{1}}\right){}^{x_{1}}x_{2}^{x_{2}}+\left(\frac{1}{x_{1}}\right){}^{x_{2}}x_{2}^{x_{2}}+\]
\[
\left(\frac{1}{x_{1}}\right){}^{x_{1}x_{2}}x_{2}^{x_{2}}+\frac{x_{2}^{x_{2}}}{x_{1}}+x_{1}x_{2}^{x_{2}}+x_{1}^{x_{1}}x_{2}^{x_{2}}+x_{1}^{x_{2}}x_{2}^{x_{2}}+x_{1}^{x_{1}x_{2}}x_{2}^{x_{2}}+x_{2}^{x_{1}x_{2}}+\left(\frac{1}{x_{1}}\right){}^{x_{1}}x_{2}^{x_{1}x_{2}}+\]
\[
\left(\frac{1}{x_{1}}\right){}^{x_{2}}x_{2}^{x_{1}x_{2}}+\left(\frac{1}{x_{1}}\right){}^{x_{1}x_{2}}x_{2}^{x_{1}x_{2}}+\frac{x_{2}^{x_{1}x_{2}}}{x_{1}}+x_{1}x_{2}^{x_{1}x_{2}}+x_{1}^{x_{1}}x_{2}^{x_{1}x_{2}}+x_{1}^{x_{2}}x_{2}^{x_{1}x_{2}}+x_{1}^{x_{1}x_{2}}x_{2}^{x_{1}x_{2}}\]

\noindent \begin{equation}
\end{equation}
It should be noted that in the resulting expression the towers associated
with rational appear only once in the sequence.

\subsection{Addition and Subtraction Algorithms.}

We assume that numbers are given in their prime tower representation
and consider the recursive algorithm expressed by 

\begin{equation}
\begin{cases}
\begin{array}{c}
n+m=\left(\sqrt{n}+\sqrt{m}\right)^{2}-2\sqrt{n}\cdot\sqrt{m}\\
n-m=\left(\sqrt{n}-\sqrt{m}\right)\cdot\left(\sqrt{n}+\sqrt{m}\right)\end{array}\end{cases}\end{equation}
 An alternative option for adding towers include sieving and checking
or solving the optimization proble induced by the identity

\begin{equation}
d\left(T_{\boldsymbol{X}}\left(m\right),\, T_{\boldsymbol{X}}\left(p\right)\right)=d\left(0,\, T_{\boldsymbol{X}}\left(n\right)\right)\Leftrightarrow T_{\boldsymbol{X}}\left(m\right)+T_{\boldsymbol{X}}\left(n\right)=T_{\boldsymbol{X}}\left(p\right)\end{equation}
 assuming that we have recovered the ordering for sufficiently many
towers we may check for candidate towers $T_{\boldsymbol{X}}\left(p\right)$
for which the following conditions are fulfilled\begin{equation}
\begin{array}{c}
T_{\boldsymbol{X}}\left(p\right)\ge T_{\boldsymbol{X}}\left(n\right)\\
and\\
T_{\boldsymbol{X}}\left(p\right)\ge T_{\boldsymbol{X}}\left(m\right)\end{array}\end{equation}
\begin{equation}
\left(T_{\boldsymbol{X}}\left(m\right)+T_{\boldsymbol{X}}\left(n\right)=T_{\boldsymbol{X}}\left(p\right)\right)\Leftrightarrow\begin{cases}
\begin{array}{c}
d\left(T_{\boldsymbol{X}}\left(m\right),\, T_{\boldsymbol{X}}\left(p\right)\right)\ge d\left(0,\, T_{\boldsymbol{X}}\left(n\right)\right)\\
d\left(T_{\boldsymbol{X}}\left(m\right),\, T_{\boldsymbol{X}}\left(p\right)\right)\le d\left(0,\, T_{\boldsymbol{X}}\left(n\right)\right)\end{array}\end{cases}\end{equation}

\subsection{Functional relation of tower progression}

Tower progression are defined by the following expression 

\begin{equation}
S_{n}\left(x\right)=1+x+x^{x}+x^{\left(x^{x}\right)}+\cdots+\left(\left.x^{\left(x^{\left(.^{.^{.\left(x^{x}\right)}}\right)}\right)}\right\} height=n\right)\end{equation}
from which it follows that 

\begin{equation}
\begin{cases}
\begin{array}{ccc}
S_{n+1}\left(x\right) & = & 1+\pmb{\mathfrak{R}}\left(x,\, S_{n}\left(x\right)\right)\\
S_{n+1}\left(x\right) & = & S_{n}\left(x\right)+\left(\left.x^{\left(x^{\left(.^{.^{.\left(x^{x}\right)}}\right)}\right)}\right\} height=n\right)\end{array}\end{cases}\end{equation}
So that tower progression is determined by the following equation

\begin{equation}
1+\pmb{\mathfrak{R}}\left(x,\, S_{n}\left(x\right)\right)-S_{n}\left(x\right)-\left(\left.x^{\left(x^{\left(.^{.^{.\left(x^{x}\right)}}\right)}\right)}\right\} height=n\right)=0\end{equation}

\section{Conclusion}

In part I we have proposed an inherently combinatorial approach to
investigating properties of numbers.

\subsection*{Acknowledgment:}

I am deeply grateful to Professor Doron Zeilberger for his insightful
comments, suggestions and encouragements. I am  grateful to Eric Rowland
who patiently and diligently thought me everything I know about Mathematica.
I am also grateful to Pavel Kuksa for insightful discussions.

\end{document}